\documentclass[12pt,leqno,]{article}
\usepackage{amsmath,amsthm,amssymb}
\usepackage{graphicx}
\usepackage{float}
\usepackage{color}
\usepackage[utf8]{inputenc}
\usepackage[hang,flushmargin]{footmisc} 
\usepackage{epsfig}
\usepackage{marvosym}

\newtheorem{theorem}{THEOREM}[section] 
\newtheorem{lemma}{Lemma}[section]

\newtheorem{remark}{REMARK}[section]

\newtheorem{propos}{Proposition}[section]

\newtheorem{examp}{Example}[section]

\newcommand\blfootnote[1]{%
  \begingroup
  \renewcommand\thefootnote{}\footnote{#1}%
  \addtocounter{footnote}{-1}%
  \endgroup
}

\numberwithin{equation}{section}

\newcommand{\dx}{\, \mathrm{d}x}

\newcommand{\ds}{\, \mathrm{d}s}
\newcommand{\dr}{\, \mathrm{d}r}
\newcommand{\dt}{\, \mathrm{d}t}

\newcommand{\rz}{\mathbb{R}^{2}}

\newcommand{\R}{\mathbb{R}}
\newcommand{\loc}{\mathrm{loc}}

\newcommand{\intom}{\int_\Omega}

\def\IntQ2Rx{\int_{{Q_{2R} (x_0)}}}

\def\O{\Omega}
\newcommand{\spt}{\operatorname{spt}}

\newcommand{\newmathop}[2]{\newcommand{#1}{\mathop{\mit{#2}}}}
\newmathop{\Minttext}{\int\limits\hspace{ -9.2em}-}
\newmathop{\Mint}{\int\hspace{ -1.0em}-}

\addtocounter{section}{0}

\title{\bf A remark on the denoising of greyscale images using energy densities with varying growth rates}
\author{M.~Fuchs \and J.~M\"uller}
\date{}

\begin{document}

\parindent2ex

\maketitle
\begin{minipage}{16cm}
AMS-classification: 49N60, 49J45, 49Q20

\end{minipage}

\begin{abstract} We prove the solvability in Sobolev spaces for a class of variational problems related to the TV-model proposed by Rudin, Osher and Fatemi in \cite{ROF} for the denoising of greyscale images. In contrast to their approach we discuss energy densities with variable growth rates depending on $|\nabla u|$ in a rather general form including functionals of $(1, p)$-growth. 
\end{abstract}

\section{Introduction}
In 1992 Rudin, Osher and Fatemi proposed (compare \cite{ROF}) to study the variational problem
\begin{equation}
\label{G1}
I_1 [w] := \intom |\nabla w| \dx + \frac{\lambda}{2} \intom |f - w|^2 \dx \to \min 
\end{equation}
as a model for the restoration of a noisy greyscale image $f$. In this setting (and throughout our paper) $\O$ is a bounded Lipschitz region in $\rz$, the function $f : \O \to \R$ represents the noisy data, for which we assume
\begin{equation}
\label{G2}
0 \le f \le 1\;\text{ a.e. on } \ \O,
\end{equation}
and $\lambda > 0$ denotes a parameter being under our disposal. As a matter of fact, problem (\ref{G1}) has to be discussed in the space $BV(\O)$ of functions with finite total variation (see, e.g., \cite{Giu} or \cite{AFP} for a definition and further properties of this class) admitting a unique solution $u$ which in addition satisfies (\ref{G2}). From the analytical point of view, the functional $I_1$ from (\ref{G1}) does not behave very nicely: the energy density $|\nabla w|$ is neither differentiable nor strictly convex (``elliptic'') so that no additional information on the minimizer $u$ are available. One common alternative used in the variational approach towards the denoising of images is to replace (\ref{G1}) by
\begin{equation}
\label{G3}
I_p [w] := \intom |\nabla w|^p \dx + \frac{\lambda}{2} \intom |w - f|^2 \dx \to \min
\end{equation}
for some power $p > 1$, where the choice $p = 2$ already occurs in the work of Arsenin and Tikhonov \cite{TA}, we refer to the monograph \cite{We} for more information on the subject including references. The natural space for problem (\ref{G3}) is the Sobolev class $W^{1,p} (\O)$ (compare \cite{Ad} for details), and from nowadays standard results on nonlinear elliptic equations (see the references stated in Chapter 3.2 of \cite{MZ}) going back to e.g. Uralt'seva, Uhlenbeck, Evans, Di Benedetto and many other prominent authors it follows that the unique solution of problem (\ref{G3}) is at least of class $C^1$ on the interior of the domain $\O$. However, from the point of view of applications, a high degree of regularity of the minimizer is not always favourable (``effect of oversmoothing''), which means that in certain cases one should discuss a linear growth model but with better ellipticity properties in comparison to the functional $I_1$. This is the subject of the papers \cite{BF1,BF2,BF3}, in which we studied the problem
\begin{equation}
\label{G4}
J_\mu [w] := \intom F_\mu (\nabla w) \dx  + \frac{\lambda}{2} \intom |w - f|^2 \dx \to \min
\end{equation}
(including even inpainting) with density
\begin{equation}
\label{G5}
F_\mu (\xi) := \Phi_\mu (|\xi|), \ \xi \in \rz,
\end{equation}
the function $\Phi_\mu : [0, \infty) \to [0, \infty)$ being defined through
\begin{equation}
\label{G6}
\Phi_\mu (t) := \int_0^t \int_0^s (1 + r)^{-\mu}\, \mathrm{d}r\, \mathrm{d}s, \;t \ge 0,
\end{equation}
with explicit formula
\begin{align}
\label{G7}
\left\{
\begin{aligned}
&\Phi_\mu (t) = \frac{1}{\mu - 1} t + \frac{1}{\mu - 1} \frac{1}{\mu - 2} (t + 1)^{- \mu + 2} - \frac{1}{\mu -  1} \frac{1}{\mu - 2} , \ \mu \not= 2 ,\\
&\Phi_2 (t) = t - \ln (1 + t), \ t \ge 0.
\end{aligned}
\right.
\end{align}
In the case $\mu > 1$ the density $F_\mu$ is of linear growth in the sense that
\begin{equation}
\label{G8}
c_1 \big(|\xi| - 1\big) \le F_\mu (\xi) \le c_2 \big(|\xi| + 1\big)\, , \ \xi \in \rz\, ,
\end{equation}
with constants $c_1, c_2 > 0$. Formally we can also consider values $\mu < 1$, but then (\ref{G4}) reduces to (\ref{G3}) for the choice $p = 2 - \mu$. The density $F_\mu$ is of class $C^2$ satisfying in case $\mu > 1$ the condition of $\mu$-ellipticity, i.e. 
\begin{equation}
\label{G9}
c_3 \big(1 + |\xi|\big)^{- \mu} |\eta|^2 \le D^2 F_\mu (\xi) (\eta, \eta) \le c_4 \big(1 + |\xi|\big)^{-1} |\eta|^2
\end{equation}
with $c_3, c_4 > 0$ and for all $\xi, \eta \in \rz$. From (\ref{G7}) it follows
\begin{equation}
\label{G10}
\lim_{\mu \to \infty} \left(\mu - 1\right) F_\mu \left(\xi\right) = \left|\xi\right|\, , \ \xi \in \rz\, ,
\end{equation}
and (\ref{G9})  together with (\ref{G10}) shows that ``$(1 - \mu) F_\mu (\nabla w)$'' is a reasonable approximation of the TV-density ``$|\nabla w|$'' occurring in problem (\ref{G1}). Moreover, it turns out that the degree of regularity of the solution $u_\mu \in \mbox{BV} (\O)$ of problem (\ref{G4}) can be controlled in terms of the parameter $\mu$. Precisely it holds 
\begin{theorem}
Let $f$ satisfy (\ref{G2}), fix $\mu>1$ and define $F_\mu$ according to (\ref{G5}), (\ref{G6}).
\begin{enumerate}
\item[a)] If $\mu < 2$, then the solution $u_\mu$ of (\ref{G4}) belongs to the Sobolev space $W^{1,1} (\O)$ and is of class $C^1$ in the interior of $\O$.
\item[b)] In case $\mu > 2$ there are simple examples of data $f$ for which $u_\mu \not\in W^{1,1} (\O)$.
\end{enumerate}
\end{theorem}
\noindent For part a) we refer to \cite{BF1,BF2,BF3,BFT,BFW2}, a discussion of b) even for the one-dimensional case $\O = (0, 1)$ can be found in \cite{FMT}.
Up to now, all our energy functionals are of uniform power growth in the sense that the regularizing part involving $\nabla w$ can be estimated from above and below by the quantity $\intom |\nabla w|^q \dx$ for some power $q \in [1, \infty)$, and the purpose of the present paper is to introduce - at least to some extend - energy functionals and densities $F$, which allow some flexibility of the growth rate, which means that the growth rate of $F (\nabla w)$ can be prescribed in terms of $|\nabla w|$. To be precise, we consider a density $F : \rz \to [0, \infty)$ of class $C^2$ satisfying $F (0) = 0$ and $D F (0) = 0$. For numbers $c_5, c_6 > 0$ and for exponents
\begin{equation}
\label{G11}
p, \mu \in (1, \infty)
\end{equation}
we assume the validity of $(\eta, \xi \in \rz)$
\begin{equation}
\label{G12}
c_5 \left(1 + |\xi|\right)^{- \mu} |\eta|^2 \le D^2 F (\xi) (\eta, \eta) \le c_6 \left(1 + |\xi|\right)^{p - 2} |\eta |^2\, ,
\end{equation}
and in Lemma 2.1 we will show that (\ref{G12}) yields the growth estimate $(c_7, \tilde{c}_7, c_8 > 0)$
\begin{equation}
\label{G13}
c_7 |\xi| - \tilde{c}_7 \le F (\xi) \le c_8 \left(|\xi|^p + 1\right)\, .
\end{equation}
The reader should note that (\ref{G12}) implies (\ref{G9}), if we allow the choice $p = 1$. An example of a density $F$ with (\ref{G12}) is given by $(\varepsilon > 0)$
\begin{equation}
\label{G14}
F (\xi) := \int^{|\xi|}_0 \int^s_0 (\varepsilon + r)^{\varphi (r) - 2} \dr \ds\, , \ \xi \in \rz\, ,
\end{equation}
for a continuous and decreasing function
\begin{align*}
\varphi : [0, \infty) \to [2 - \mu, p], \quad \varphi (0) = p,\quad \lim\limits_{r \to \infty} \varphi (r) = 2 - \mu.
\end{align*}
A discussion of (\ref{G14}) together with further examples can be found in Section 5. Assuming (\ref{G12}) we then look at the variational problem
\begin{equation}
\label{G15}
J [w] := \intom F (\nabla w) \dx + \frac{\lambda}{2} \int_{ \O - D} |w - f|^2 \dx\, ,
\end{equation}
where $D$ is a measurable subset of $\O$ such that
\begin{equation}
\label{G16}
0 \le  {\cal{L}}^2 (D) < {\cal{L}}^2 (\O)\, ,
\end{equation}
i.e. we study an inpainting problem combined with simultaneous denoising, where $D$ is the inpainting region and the choice $D = \emptyset$ corresponds to the case of pure denoising.
We have the following results:
\begin{theorem}
Let (\ref{G2}), (\ref{G11}) and (\ref{G12}) hold together with (\ref{G16}). Assume in addition that
\begin{equation}
\label{G17}
\mu, p < 2\, .
\end{equation}
Then the variational problem
\begin{equation}
\label{G18}
J [w] \to \min \mbox{in} \ W^{1,1} (\O)
\end{equation}
with $J$ defined in (\ref{G15}) admits a unique solution $u$. This solution additionally satisfies $0 \le u \le 1$ a.e. on $\O$ as well as $u \in W^{1,s}_\loc (\O)$ for any finite $s$.
\end{theorem}
\begin{remark}
Once having established the local higher integrability result $|\nabla u| \in L^s_\loc (\O), s < \infty$, we think that actually $u \in C^{1, \alpha} (\O)$, $0 < \alpha < 1$, can be deduced along similar lines as in \cite{BFT}, where densities $F$ satisfying (\ref{G9}) for some exponent $\mu \in (1,2)$ are considered.
\end{remark}
\begin{remark} Energy densities $F$, for which
\begin{equation}
\label{G19}
c_9 |\nabla w|^s - \tilde{c}_9 \le F (\nabla w) \le c_{10} \left(\left|\nabla w\right|^q + 1\right)
\end{equation}
 holds or for which an appropriate variant of (\ref{G12}) is true, have been extensively discussed for instance in the papers \cite{Ma2,Ma3,Ma4,Ma6,ELM1,ELM2,BF10,BF11,BF12} dealing even with the higher-dimensional case including vector-valued functions. Roughly speaking it is shown in the above mentioned papers and the references quoted therein, that (\ref{G19}) provides some additional regularity of (local) minimizers, provided $s > 1$ and $q$ is not too far away from $s$, we refer to \cite{Bi} for a survey. Recalling that (\ref{G12}) implies (\ref{G13}), Theorem 1.2 covers the case ``$s = 1$'', and (\ref{G17}) expresses the fact that the upper bound $p$ satisfies ``$p < 2 s$''. Note that the latter requirement turns out to be a sufficient condition for the regularity of minimizers in the setting of \cite{BF11}.
\end{remark}
\begin{remark}
Variational problems of mixed linear/superlinear growth are the subject of Section 6 in \cite{Bi}. Here the density $F$ is of splitting form in the sense that
\begin{equation}
\label{G20}
F (\nabla w) = F (\partial_1 w \  \partial_2 w) = F_1 (\partial_1 w) + F_2 (\partial_2 w)
\end{equation}
with $F_1$ growing linearly in $|\partial_1 w|$, whereas $F_2 (\partial_2 w)$ behaves as $|\partial_2 w|^p$ with power $p > 1$. From the point of view of image restoration condition (\ref{G20}) seems to be unnatural, however, if $F_1$ satisfies (\ref{G9}) with $\mu \in (1,2)$ and if $p < 2$, then regularity results are available, thus our hypothesis (\ref{G17}) naturally occurs in the splitting case (\ref{G20}).
\end{remark}

Next let $\rho : [0, \infty) \to [0, \infty)$ denote a function of class $C^1$ being strictly increasing and strictly convex, e.g. $\rho (t) = \sqrt{1+ t^2} - 1$, and let
\begin{equation}
\label{G21} K [w] := \intom F (\nabla w) \dx + \int_{\O - D} \rho \left(|w - f|\right) \dx\, ,
\end{equation}
which means that we consider more general data terms.
\begin{theorem}
With $\rho$ from above let $f$, $F$ and $D$ satisfy (\ref{G2}), (\ref{G12}) and (\ref{G16}), respectively, and assume in addition that $\displaystyle\limsup_{t \to \infty} \frac{\rho (t)}{t^m} < \infty$ for some $m \ge 1$. Moreover, let
\begin{align}
\label{G22}
& 1 < \mu < 3/2\, ,\\
\label{G23}
&1 < p < \mu\, .
\end{align}
Then the variational problem
\begin{equation}
\label{G24}
K [w ] \to \min \text{ in } \ W^{1,1} (\O)
\end{equation}
with $K$ from (\ref{G21}) has a unique solution $u$. It holds $0 \le u \le 1$ a.e. on $\O$, moreover, $|\nabla u|$ is in $L^s_\loc (\O)$ for any finite $s$. If the density $F$ is balanced in the sense that 
\begin{equation}
\label{G25} 
\left|D^2 F (\xi)\right| \left|\xi\right|^2 \le c_{11} \left(F\left(\xi\right) + 1\right), \ \xi \in \rz\, ,
\end{equation}
holds for some constant, then (\ref{G23}) can be replaced by the requirement $p \in (1,2)$ (compare (\ref{G17})).
\end{theorem}
\begin{remark}
We conjecture that in the balanced case (\ref{G25}) the results of Theorem 1.2 and 1.3 extend to any exponent $p \ge 2$, we refer to Remark 3.1.
\end{remark}

Our paper is organized as follows: in Section 2 we collect some preliminary material and discuss regularized problems approximating (\ref{G18}) and (\ref{G24}). Section 3 is devoted to the proof of Theorem 1.2, and Theorem 1.3 is established in Section 4. Finally, in Section 5 we present some examples of densities $F$ satisfying (\ref{G12}) including the model from (\ref{G14}).

\section{Some preliminary results and discussion of regularized problems}
We start with a growth estimate for densities $F$ satisfying (\ref{G12}). 
\begin{lemma}
Suppose that we have the ellipticity condition (\ref{G12}) for $F : \rz \to [0, \infty)$ with exponents $p, \mu$ according to (\ref{G11}). Then $F$ is of $(1, p)$-growth in the sense of inequality (\ref{G13}).
\end{lemma}
\noindent\textit{Proof.} We just consider the case $p \ge 2$. For $ p < 2$ the following arguments can be easily adjusted. We recall that $F$ should satisfy $F (0) = 0$, $D F (0) = 0$, thus we obtain from Taylor's theorem (applied to $t \mapsto F (t \xi)$)
\begin{equation}
\label{H1} 
F (\xi) = \int^1_0 (1 - t) D^2 F (t\xi)(\xi, \xi) \dt, \ \xi \in \rz\, .
\end{equation}
Applying (\ref{G12}) to the r.h.s. of (\ref{H1}) we find
\begin{equation}
\label{H2}
c_{12} \int^1_0 (1 - t) (1 + t |\xi|)^{- \mu} \dt |\xi|^2 \le F (\xi) \le c_{13} \int^1_0 (1 - t) (1 + t |\xi|)^{p - 2} \dt |\xi|^2\, ,
\end{equation}
and from $(1 + t |\xi|)^{p - 2} |\xi|^2 \le (1 + |\xi|)^p$ (in case $p \ge 2$) we immediately deduce the second inequality in (\ref{G13}). If $|\xi| \le 2$, then the first inequality in (\ref{G13}) is obvious by an appropriate choice of $c_7, \tilde{c}_7 > 0$. In case $|\xi| \ge 2$ we observe for the l.h.s. of (\ref{H2})
\begin{align*}
&c_{12} \displaystyle{\int^1_0}(1 - t) (1 + t |\xi|)^{- \mu} \dt |\xi|^2 \ge c_{12} \displaystyle{ \int^{1/|\xi|}_0} (1 - t) (1 + t |\xi|)^{- \mu} \dt |\xi|^2\\
& \ge c_{12} \int^{1/|\xi|}_0 (1 - t) (1+ 1)^{- \mu} \dt |\xi|^2\ge c_{14}\displaystyle{\int^{1/|\xi|}_{1/2|\xi|}} (1 - t) \dt |\xi|^2  \\
&\ge c_{14} \int^{1/|\xi|}_{1/2|\xi|} \left(1 - \frac{1}{|\xi|}\right) \dt |\xi|^2 \ge c_{14}  \displaystyle{\int^{1/|\xi|}_{1/2|\xi|}} \frac{1}{2} \dt |\xi|^2 = c_{15} |\xi|, \\
\end{align*}
thus the first inequality of (\ref{G13}) extends to the case $|\xi|\ge 2$ after adjusting $c_7, \tilde{c}_7$. \qed

\begin{remark}
The requirement $DF (0) = 0$ is essential for deducing the lower bound on $F$ stated in (\ref{G13}) from the condition of $\mu$-ellipticity, i.e. from the first inequality in (\ref{G12}).
\end{remark}
\begin{lemma}
Under the conditions on the data stated in Theorem 1.2 and 1.3, respectively, but for arbitrary choices of $p , \mu \in (1, \infty)$, the variational problems (\ref{G18}) and (\ref{G24}) admit at most one solution $u \in W^{1,1,} (\O)$. We have
\begin{equation}
\label{H3}
0 \le u \le 1 \ \ \mbox{a.e.  on} \ \O\, .
\end{equation}
\end{lemma}
\noindent \textit{Proof.}  From ``strict convexity'' (for the density $F$ this property follows from the first inequality in (\ref{G12})) we get
\[
\left\{
\begin{array}{lll}
\nabla u = \nabla v &\mbox{a.e.} & \mbox{on} \ \O\, ,\\
u = v &\mbox{a.e.}& \mbox{on} \ \O - D
\end{array}
\right.
\]
for minimizers $u, v \in W^{1,1} (\O)$. But then $u = v$ is a consequence of (\ref{G16}). Replacing $u$ by $\min (u, 1)$ and $\max (u, 0)$ we see by an elementary calculation (compare, e.g., \cite{BF2}) that (\ref{H3}) holds for the minimizer $u$, since otherwise we could decrease the energy. \qed

\noindent During the proofs of Theorem 1.2 and 1.3 we will essentially benefit from
\begin{lemma}
Suppose that we are in the situation of Theorem 1.2 or 1.3, where  here we allow in both cases exponents $p \in (1,2)$ and $\mu \in (1, \infty)$. For $\delta > 0$ let $u_\delta \in W^{1,2} (\O)$ denote the solution of either 
\begin{equation}
\label{G18D}\tag*{$(1.18)_\delta$}
J_\delta [w] := \frac{\delta}{2} \intom \vert\nabla w\vert^2 \dx + J [w] \to \min \ \mbox{in} \   W^{1,2} (\O) 
\end{equation}
or
\begin{equation}
\label{G24D}\tag*{$(1.24)_\delta$}
K_\delta [w] := \frac{\delta}{2} \intom \vert\nabla w\vert^2 \dx + K [w] \to \min \ \mbox{in} \ W^{1,2} (\O) 
\end{equation}
with $J$ and $K$ from (\ref{G15}) and (\ref{G21}), respectively. It holds:
\begin{enumerate}
\item[i)] $0 \le u_\delta \le 1$ \ \mbox{a.e.\  on} $\O$.
\item[ii)] The functions $u_\delta$ are of class $W^{2,2}_\loc (\O) \cap W^{1, \infty}_\loc (\O)$.
\item[iii)] We have the uniform bound $\displaystyle \sup_{\delta > 0} \|u_\delta\|_{W^{1,1}(\O)} < \infty$.
\item[iv)] Suppose that we can find an exponent $q > 1$ such that for each subdomain $\Omega^*\Subset\Omega$
\begin{equation}
\label{H4}
\sup_{\delta > 0} \int_{\O^\ast} \vert \nabla u_{\delta} \vert^q \dx \le c_{16} (\O^{\ast}) < \infty\, .
\end{equation}
Then $u_\delta \to : u$ in $L^1 (\O) \cap W^{1,q}_\loc (\O)$ as $\delta \to 0$ for a function $u \in W^{1,1} (\O)$, and $u$ solves the variational problem (\ref{G18}), respectively (\ref{G24}).
\end{enumerate}
\end{lemma}
\noindent\textit{Proof.}  i) follows as inequality (\ref{H3}) in Lemma 2.2, ii) is immediate from elliptic regularity theory, and iii) is a consequence of the first inequality in (\ref{G13}). Let us discuss iv):  from i), iii) and assumption (\ref{H4}) we deduce the existence of $u \in \mbox{BV} (\O) \cap W^{1,q}_\loc (\O) \subset W^{1,1} (\O)$ such that
\begin{align}
&\label{H5}  u_\delta \to u  \text{ in } \ L^1 (\O)  \text{ and a.e. },\\
&\label{H6}  u_\delta \rightharpoondown u  \text{ in }  W^{1,q}_\loc (\O)
\end{align}
(at least for a subsequence) as $\delta \to 0$. From De Giorgi's theorem on lower semicontinuity (see, e.g., \cite{Gia} Theorem 2.3, p.18) we see that (\ref{H5}) and (\ref{H6}) yield
\begin{equation}
\label{H7}
J [u] \le \liminf_{\delta \to 0} J [u_\delta]\, ,
\end{equation}
if we are in the situation of Theorem 1.2, whereas
\begin{equation}
\label{H8}
K [u] \le \liminf_{\delta \to 0} K [u_\delta]
\end{equation}
in the setting of Theorem 1.3. Since for $v \in W^{1,2} (\O)$ it holds
\[
J_\delta [u_\delta] \le J_\delta [v] \stackrel{\delta \to 0}{\longrightarrow} J [v]\, ,
\]
we obtain from (\ref{H7}) (recall the definition of $J_\delta$ in \ref{G18D})
\begin{equation}
\label{H9}
J [u] \le J [v]\, ,
\end{equation}
and by approximation $(W^{1,2} (\O) \ni v_k \to v$ in $W^{1,1} (\O))$, inequality (\ref{H9}) extends to $v \in W^{1,1} (\O)$. If the $u_\delta$ are the solutions of problem \ref{G24D}, then by the same arguments it follows
\begin{equation}
\label{H10}
K [u] \le K [v], \ v \in W^{1,2} (\O)\, .
\end{equation}
Consider $v \in W^{1,1} (\O)$. In case $K [v] = + \infty$, i.e.
\[
\int_{\O - D} \rho (|v - f|) \dx = + \infty\, ,
\]
there is nothing to prove. In the other case, due to the growth of $\rho$ at infinity and by (\ref{G2}), we see that $v$ is in the space $L^m (\O - D)$ and according to \cite{FT}, Lemma 2.1, we find a sequence $v_k \in C^\infty (\overline{\O})$ such that 
\[
\|v_k - v\|_{W^{1,1} (\O)} + \| v_k - v\|_{L^m (\O - D)} \longrightarrow 0
\]
as $ k \to \infty$, hence $K [v_k] \to K [v]$, and since $K [u] \le K [v_k]$ by (\ref{H10}), we finally have shown that $u$ solves (\ref{G24}). \qed

\section{Proof of Theorem 1.2}
In this section we assume that all the hypotheses of Theorem 1.2 are valid and define $u_\delta$ as in Lemma 2.3 as the unique solution of problem \ref {G18D}. Let $F_\delta (\xi) := \frac{\delta}{2} |\xi|^2 + F (\xi), \ \xi \in \rz$. For $\eta \in C^1_0 (\O)$ with $0 \le \eta \le 1$ we have (by passing to the differentiated version of the Euler equation associated to \ref{G18D} and by quoting Lemma 2.3 ii))
\begin{equation}
\label{I1}
\intom D^2 F_\delta\left(\nabla u_\delta\right)\left(\partial_\alpha \nabla u_\delta, \nabla [\eta^2 \partial_\alpha u_\delta]\right) \dx = \lambda \int_{\O - D} (u _\delta - f) \partial_\alpha \left(\eta^2 \partial_\alpha u_\delta\right) \dx \, ,
\end{equation}
where here and it what follows the sum in taken w.r.t. $\alpha = 1,2$. It holds
\begin{align*}
\text{r.h.s. of (\ref{I1})}=
 \lambda  \intom u_\delta \partial_\alpha \left(\eta^2 \partial_\alpha u_\delta\right) \dx- \lambda \int_D u_\delta \partial_\alpha \left(\eta^2 \partial_\alpha u_\delta\right) \dx \\
 - \lambda \int_{\O - D} f \partial_\alpha \left(\eta^2 \partial_\alpha u_\delta \right) \dx =: T_1 - T_2 - T_3,
\end{align*}
\begin{gather*}
T_1 = - \lambda \intom \eta^2 |\nabla u_\delta|^2 \dx,\\
|T_2| + |T_3| \le c_{17} \left\{\intom \eta |\nabla \eta| |\nabla u_\delta| \dx + \intom \eta^2 |\nabla^2 u_\delta| \dx \right\},
\end{gather*}
where we have used (\ref{G2}) as well as Lemma 2.3 i), $c_k$ denoting a positive constant independent of $\delta$. Recalling in addition Lemma 2.3 iii) we get from (\ref{I1})
\begin{eqnarray}
\label{I2}
\lefteqn{\intom D^2 F_\delta (\nabla u_\delta) \left(\partial_\alpha \nabla u_\delta, \nabla \left[\eta^2 \partial_\alpha u_\delta\right]\right) \dx + \lambda \intom \eta^2 |\nabla u_\delta|^2 \dx}\\
&& \le c_{18} \left\{\|\nabla \eta \|_{L^\infty (\O)} + \intom \eta^2 \left|\nabla^2 u_\delta\right|\dx \right\}\, .\nonumber
\end{eqnarray}
Applying the Cauchy-Schwarz inequality to the bilinear form $D^2 F_\delta (\nabla u_\delta)$ and using Young's inequality, the estimate (\ref{I2}) yields
\begin{eqnarray*}
\lefteqn{\intom \eta^2  D^2 F_\delta (\nabla u_\delta) \left(\partial_\alpha \nabla u_\delta, \partial_\alpha \nabla u_\delta\right) \dx + \intom \eta^2 |\nabla u_\delta|^2 \dx}\\
&& \le c_{19} \left\{\intom  D^2 F_\delta (\nabla u_\delta) \left(\nabla \eta, \nabla \eta\right) |\nabla u_\delta|^2 \dx+ \intom \eta^2 \left|\nabla^2 u_\delta\right| \dx + \| \nabla \eta \|_{L^\infty (\O)} \right\}\, ,
\end{eqnarray*}
hence using (\ref{G12}) for $D^2 F$ (dropping the $\delta$-term on the l.h.s.)
\begin{eqnarray*}
\lefteqn{\intom \eta^2 \left(1 + |\nabla u_\delta|\right)^{- \mu} \left|\nabla^2 u_\delta\right|^2 \dx + \intom \eta^2 |\nabla u_\delta|^2 \dx}\\
&& \le c_{20} \left\{\|\nabla \eta\|^2_{L^\infty (\O)} \delta \intom \left|\nabla u_\delta\right|^2 \dx +\left\|\nabla \eta \right\|^2_{L^\infty (\O)} \int_{\spt \eta} \left(1 + |\nabla u_\delta|\right)^p \dx\right.\\
&&\left. + \left\| \nabla \eta \right\|_{L^\infty (\O)} + \intom \eta^2 \left|\nabla^2 u_\delta\right| \dx \right\}\, .
\end{eqnarray*}
We remark the validity of $\displaystyle\sup_{\delta > 0} \,\delta \intom \left|\nabla u_\delta\right|^2 \dx < \infty$ and 
assume w.l.g. $\|\nabla \eta\|_{L^\infty (\O)} \ge 1$. Then we obtain
\begin{eqnarray}
\label{I3}
\hspace*{1cm}\lefteqn{\intom \eta^2 \left(1 + |\nabla u_\delta|\right)^{- \mu} \left|\nabla^2 u_\delta\right|^2 \dx + \intom \eta^2 |\nabla u_\delta|^2 \dx}\\
&& \le c_{21} \left\{\|\nabla \eta\|^2_{L^\infty (\O)} \int_{\spt \eta} \left(1 + |\nabla u_\delta|\right)^p \dx + \intom \eta^2 \left|\nabla^2 u_\delta\right| \dx + \left\| \nabla \eta \right\|^2_{L^\infty (\O)} \right\}\, .\nonumber
\end{eqnarray}
On the r.h.s. of (\ref{I3}) we use Young's inequality twice recalling (\ref{G17}) and (\ref{G11}):
\begin{align*}
&\left\|\nabla \eta\right\|^2_{L^\infty (\O)} \int_{\spt \eta} \left(1 + |\nabla u_\delta|\right)^p \dx\le \tau  \int_{\spt \eta} |\nabla u_\delta|^2 \dx + c_{22} (\tau) \left\|\nabla \eta \right\|^{\frac{4}{2 - p}}_{L^\infty (\O)},\\
& \intom \eta^2 \left|\nabla^2 u_\delta\right| \dx \le \varepsilon \intom \eta^2 \left(1 + |\nabla u_\delta|\right)^{- \mu} \left|\nabla^2 u_\delta\right|^2 \dx + c_{23} (\varepsilon) \intom \left(1 + |\nabla u_\delta|\right)^\mu  \eta^2 \dx\\
& \le \varepsilon \intom \eta^2 \left(1 + |\nabla u_\delta|\right)^{- \mu} \left|\nabla^2 u_\delta\right|^2 \dx + \varepsilon \intom \eta^2 |\nabla u_\delta|^2 \dx  +  c_{24} (\varepsilon).
\end{align*}
Inserting these estimates into (\ref{I3}), choosing $\eta$ such that $\eta \equiv 1$ on $B_{r_1} (x_0), \eta \equiv 0$ outside $B_{r_2} (x_0), B_{r_1} (x_0) \subset B_{r_2} (x_0) \Subset \O$, we obtain after appropriate choice of $\varepsilon$ and $\tau$
\begin{equation}
\label{I4}
\int_{B_{r_1} (x_0)}  |\nabla u_\delta|^2 \dx \le \frac{1}{2} \int_{B_{r_2} (x_0)} |\nabla u_\delta|^2 \dx + c_{25} \left(\left(r_2 - r_1\right)^{- \alpha} + 1\right)\, ,
\end{equation}
where for the moment we just neglect $\intom \eta^2 (1 + |\nabla u_\delta|)^{- \mu} |\nabla^2 u_\delta|^2 \dx$ and $\alpha$ denotes a suitable positive number. Applying Lemma 3.1, p.161, from \cite{Gia} to estimate (\ref{I4}) we find that (\ref{H4}) from Lemma 2.3 holds with the choice $q = 2$, and we can quote iv) of Lemma 2.3 yielding a unique $W^{1,1} (\O)$- solution $u$ of (\ref{G18}).\\
Going back to (\ref{I3}), recalling the estimates stated after (\ref{I3}) and applying our bound (\ref{H4}) valid for $q = 2$, it follows
\[
\int_{\O^*} \left|\nabla^2 u_\delta\right|^2 \left(1 + |\nabla u_\delta|\right)^{- \mu} \dx \le c_{26} (\O^\ast) < \infty
\]
for any $\O^\ast \Subset \O$, thus ($\varphi_\delta := (1 + |\nabla u_\delta|)^{1 - \mu/2}$)
\[
\|\varphi_\delta\|_{W^{1,2} (\O^*)} \le c_{27} (\O^*) < \infty\, ,
\]
which by Sobolev's theorem implies
\begin{equation}
\label{I5}
\|\nabla u_\delta\|_{L^s (\O^*)} \le c_{28} (s, \O^*)
\end{equation}
for any $s < \infty$. This proves the last claim of Theorem 1.2. \qed

\begin{remark}
Suppose that $F$ satisfies the condition (\ref{G25}). In this case we estimate 
\[
D^2 F_\delta (\nabla u_\delta) (\nabla \eta, \nabla \eta) \left|\nabla u_\delta\right|^2 \le c_{29} \left(F_\delta (\nabla u_\delta) + 1\right)
\]
and observe $\displaystyle\sup_{\delta > 0} \intom F_\delta (\nabla u_\delta) \dx < \infty$. Thus we can replace $\int_{\spt \eta} \left(1 + |\nabla u_\delta|\right)^p \dx$ in (\ref{I3}) through a constant ending up with
\[
\intom \eta^2 \left(1 + |\nabla u_\delta|\right)^{- \mu} \left|\nabla^2 u_\delta\right|^2 \dx + \intom \eta^2 \left|\nabla u_\delta\right|^2 \dx \le c_{30} (\eta)\, ,
\]
hence we obtain (\ref{I5}) just assuming $\mu \in (1,2)$. Thus the bound (\ref{G17}) imposed on $p$ at this stage does not enter, however during our proof we work with the quadratic regularization \ref{G18D}, which requires $p \le 2$. In other words: under the assumption (\ref{G25}) the claims of Theorem 1.2 extend to exponents $p > 2$ (keeping the bound $1 < \mu < 2$) and a proof can be carried out by working with the regularization
\[
\delta \intom \left(1 + |\nabla w|^2\right)^{\overline{p}/2} \dx + J [w] \to \min \ \mbox{in} \ W^{1, \overline{p}} (\O)
\]
for some exponent $\overline{p} > p$. We leave the details to the reader. 
\end{remark}

\section{Proof of Theorem 1.3}
Let the assumptions of Theorem 1.3 hold. In place of equation (\ref{I1}) we have
\begin{eqnarray}
\label{J1}
\lefteqn{\intom D^2 F_\delta (\nabla u_\delta) \left(\partial_\alpha \nabla u_\delta, \nabla \left[\eta^2 \partial_\alpha u_\delta \right]\right)\dx}\\
&& = \int_{\O - D} \rho' \left(|u_\delta - f|\right) \frac{u_\delta - f}{|u_\delta - f|} \partial_\alpha \left(\eta^2 \partial_\alpha u_\delta\right) \dx\, ,\nonumber
\end{eqnarray}
where $u_\delta$ is the solution of problem \ref{G24D} (see Lemma 2.3). From (\ref{G2}) and Lemma 2.3 i) it follows
\[
\mbox{r.h.s. of (\ref{J1})} \ \le c_{31} \intom \left|\partial_\alpha \left(\eta^2 \partial_\alpha u_\delta \right)\right| \dx
\]
and clearly (recall Lemma 2.3 iii))
\begin{equation}
\label{J2}
\intom \left|\partial_\alpha \left(\eta^2 \partial_\alpha u_\delta\right)\right| \dx \le c_{32} (\eta) + c_{33} \intom \eta^2 \left|\nabla^2 u_\delta\right| \dx\, ,
\end{equation}
where we use the symbol $c_k (\eta)$ to denote constants proportional to $\|\nabla \eta\|^\alpha_{L^\infty (\O)}$ for some positive exponent $\alpha$. Applying Young's inequality to the  integral on the r.h.s. of (\ref{J2}) and discussing the l.h.s. of (\ref{J1}) as done after (\ref{I2}) we find
\begin{eqnarray}
\label{J3}
\lefteqn{\intom \eta^2 \left(1 + |\nabla u_\delta|\right)^{- \mu} \left|\nabla^2 u_\delta\right|^2 \dx}\\
&& \le c_{34} \intom \eta^2 \left(1 + |\nabla u_\delta|\right)^\mu \dx + c_{35} (\eta) \int_{\spt \eta} \left(1 + |\nabla u_\delta|\right)^p \dx \, .\nonumber
\end{eqnarray}
We specify $\eta$ as in Section 3 and let
\[
\varphi_\delta := \left(1 + |\nabla u_\delta|\right)^{1 - \mu/2}, \ \Psi_\delta := \left(1 + |\nabla u_\delta|\right)^{\mu/2}\, .
\]
Then (\ref{J3}) shows (with suitable $\alpha_1 > 0$)
\begin{align}
\label{J4}
\intom \eta^2 |\nabla \varphi_\delta|^2 \dx \le c_{36} \left\{\intom \eta^2 \Psi^2_\delta \dx + (r_2 - r_1)^{-\alpha_1} \int_{B_{r_2}(x_0)} \left(1 + |\nabla u_\delta|\right)^p \dx \right\}.
\end{align}

\noindent Next we observe (quoting Sobolev's inequality)
\begin{align}
\label{J5}
\begin{split}
\intom (\eta \Psi)^2 \dx \le c_{37} \left(\intom \left|\nabla (\eta \Psi_\delta)\right| \dx \right)^2\le c_{38} \left[\intom |\nabla \eta| \Psi_\delta \dx + \intom \eta \left|\nabla \Psi_\delta \right| \dx \right]^2 \\
 \le c_{39} (\nabla \eta) + c_{40} \left(\intom |\nabla \Psi_\delta|\eta \dx \right)^2,
 \end{split}
\end{align}
where we have used that $\sup_{\delta > 0} \intom \Psi_\delta \dx < \infty$ on account of Lemma 2.3 iii). We discuss the remaining integral on the r.h.s. of (\ref{J5}) observing that $\Psi_\delta = \varphi^{\mu/(2-\mu)}_\delta$ and using H\"older's inequality:
\begin{align*}
&\intom \eta |\nabla \Psi_\delta| \dx \le c_{41} \intom \eta|\nabla \varphi_\delta| \varphi_\delta^{\frac{\mu}{2 - \mu} - 1} \dx\\
&\le c_{42} \left(\intom \eta^2 \left|\nabla \varphi_\delta \right|^2 \dx \right)^{1/2} 
\left(\int_{B_{r_2} (x_0)} \varphi^{2 \frac{2 \mu - 2}{2 - \mu}}_\delta \dx\right)^{1/2}\, .
\end{align*}
We have \[ \varphi^{2 \frac{2 \mu - 2}{2 - \mu}}_\delta = \left(1 + |\nabla u_\delta|\right)^{2 \mu - 2}\]
with exponent $2 \mu - 2 \in (0, 1)$, which follows from (\ref{G22}). Quoting Lemma 2.3 iii) one more time, another application of H\"older's inequality gives (for some $\alpha_2 > 0$)
\begin{equation}
\label{J6}
\intom \eta |\nabla \Psi_\delta| \dx \le c_{43}\ r^{\alpha_2}_2 \left(\intom \eta^2 \left|\nabla \varphi_\delta\right|^2 \dx \right)^{1/2}\, .
\end{equation}
We insert (\ref{J6}) into (\ref{J5}) giving the bound
\begin{equation}
\label{J7}
\intom \left(\eta\Psi_\delta\right)^2 \dx \le c_{44} \ (\nabla \eta) + c_{45} \ r^{2 \alpha_2}_2 \intom \eta^2 |\nabla \varphi_\delta|^2 \dx\, .
\end{equation}
With (\ref{J7}) we return to (\ref{J4}) and assume that the radius $r_2$ is sufficiently small, thus
\begin{equation}
\label{J8}
\intom \eta^2 |\nabla \varphi_\delta|^2 \dx \le c_{46} (r_2 - r_1)^{-\alpha_3} \left(1 + \int_{B_{r_2}(x_0)} \left(1 + |\nabla u_\delta|\right)^p \dx\right)\, .
\end{equation}
Up to now we have not used our hypothesis (\ref{G23}), which enters next:
\begin{eqnarray*}
\lefteqn{\int_{B_{r_1} (x_0)} \left(1 + |\nabla u_\delta|\right)^\mu \dx  \le \intom \left(\eta \Psi_\delta \right)^2 \dx}\\
& \stackrel{{\text{(\ref{J7}), (\ref{J8})}}}{\le}& c_{47} \left[ \left(r_2 - r_1\right)^{- \alpha_4} + \left(r_2 - r_1\right)^{- \alpha_5} \int_{B_{r_2} (x_0)} \left(1 + |\nabla u_\delta|\right)^p \dx \right]\\
&\le& c_{48} \left(r_2 - r_1\right)^{- \alpha_6} + \frac{1}{2} \int_{B_{r_2} (x_0)} \left(1 + |\nabla u_\delta|\right)^\mu \dx \, ,
\end{eqnarray*}
where in the last estimate we applied H\"older's inequality and use the smallness of $r_2$ to get the factor $1/2$. As outlined after (\ref{I4}) we deduce (\ref{H4}) with value $q := \mu$. Moreover, using this information in (\ref{J4}), we see 
\[
\sup_{\delta > 0} \|\varphi_\delta\|_{W^{1,2} (\O^*)} < \infty
\]
for any subdomain $\O^* \Subset \O$, thus (\ref{I5}) holds, and we get all the results of Theorem 1.3 as described in Section 3, where for the balancing case we refer to Remark 3.1. \qed

\section{Examples}
In this section we focus on energy densities depending on the modulus of $\nabla u$, a situation for which the following observations are helpful.
\begin{propos}
Let $g: [0, \infty) \to [0, \infty)$ denote a $C^2$-function for which $g (0) = g' (0) = 0, \ g'' \ge 0$. Then 
\begin{equation}
\label{K1}
G : \rz \to [0, \infty), \ G(\xi) := g (|\xi|)\, ,
\end{equation}
is a convex function of class $C^2$ for which $G (0) = 0 $, $DG (0) = 0$ and
\begin{align}
\label{K2}
\begin{split}
\min \left\{g'' \left(|\xi|\right), \frac{1}{|\xi|} g' \left(|\xi|\right)\right\} |\eta|^2 \le D^2 G (\xi) (\eta, \eta)\\
 \le \max \left\{g'' \left(|\xi|\right), \frac{1}{|\xi|} g' \left(|\xi|\right)\right\} |\eta|^2, \ \xi, \eta \in \rz.
\end{split}
\end{align}
\end{propos}

\noindent \textit{Proof.} We just note that (\ref{K2}) follows from the formula 
\[
D^2 G (\xi) (\eta, \eta) = \frac{1}{|\xi|} g' (|\xi|) \left[|\eta|^2 - \frac{(\eta \cdot \xi)^2}{|\xi|^2} \right] + g'' (|\xi|) \frac{(\eta \cdot \xi)^2}{|\xi|^2}\, .
\] \qed

\noindent If $G : \rz \to [0, \infty)$ is a non-negative function of class $C^2$, we recall the balancing condition (see (\ref{G25})):
\begin{equation}
\label{K3}
\left|D^2 G (\xi)\right| |\xi|^2 \le c_{49} \left(G \left(\xi\right) + 1\right), \ \xi \in \rz\, .
\end{equation}

\begin{propos}
Let $g$ satisfy the assumptions of Proposition 5.1. Assume further that
\begin{equation}
\label{K4}
t^2 \max \left\{g'' (t), \frac{1}{t} g' (t) \right\} \le c_{50} \left(g (t) + 1\right), \ t \ge 0\, .
\end{equation}
Then $G$ from (\ref{K1}) satisfies (\ref{K3}).
\end{propos}

\noindent \textit{Proof.}  This is an immediate consequence of (\ref{K2}) and (\ref{K4}). \qed

\begin{propos}
If $g$ is a function as in Proposition 5.1 such that
\begin{equation}
\label{K5}
t g'' (t) \le c_{51} \ g' (t), \ t \ge 0\, ,
\end{equation}
and if $g$ satisfies the $(\Delta 2)$-condition, i.e.
\begin{equation}
\label{K6}
 g  (2t) \le c_{52} \ g (t), \ t \ge 0\, ,
\end{equation}
then inequality (\ref{K4}) holds. Thus $G (\xi) := g (|\xi|)$, $\xi \in \rz$, is balanced in the sense of (\ref{K3}).
\end{propos}

\noindent \textit{Proof.} We have (recalling $g (0) = 0$ as well as $g' (t) \ge 0$)
\[
g (t) = \int^t_0 g' (s) \ds \ge \int^t_{t/2} g' (s) \ds \ge \frac{t}{2} g' (t/2)\, ,
\]
where the last inequality follows from the fact that $g'$ is increasing on account of our hypothesis $g'' \ge 0$.
Thus we get 
\[
t \, g' (t) \le g (2 t) \le c_{52} \ g (t)
\]
(see (\ref{K6})), and together with (\ref{K5}) we arrive at (\ref{K4}). \qed

All our energy densities ``$G (\nabla u)$'' discussed below are of the principal form (compare (\ref{G6}) and (\ref{G14}))
\begin{equation}
\label{K7}
G (\xi) := g (|\xi|) := \int^{|\xi|}_0 \int^s_0 \omega (r) \dr\ds, \ \xi \in \rz\, ,
\end{equation}
for a continuous function $\omega : [0, \infty) \to [0, \infty)$ such that $\omega (s) > 0$ for $s > 0$. Note that $g$ is strictly increasing and strictly convex implying the strict convexity of $G$ on $\rz$.

\begin{examp}
Consider a continuous function $\eta : [0, \infty) \to [0, 1]$ and define $G$ according to (\ref{K7}) with the choice
\begin{equation}
\label{K8}
\omega (t) := \eta (t) (1 + t)^{- \mu } + \left(1 - \eta (t)\right) (1 + t)^{p - 2}, \ t \ge 0\, ,
\end{equation}
for exponents $p, \mu \in (1, \infty)$.
\end{examp}
\begin{propos}
The density $G$ satisfies (\ref {G12}).
\end{propos}

\noindent\textit{Proof.} It holds on account of $\eta (t) \in [0, 1]$ and $(1+t)^{- \mu} \le (1 + t)^{p - 2}$ for any $t \ge 0$
\begin{eqnarray*}
\lefteqn{g'' (t) \stackrel{{\text{(\ref{K7})}}}{=} \omega (t) \stackrel{{\text{(\ref{K8})}}}{\le} (1 + t)^{p - 2}\, ,}\\[2ex]
&&\hspace*{-0.75cm}\displaystyle{\frac{1}{t}} g' (t) \stackrel{{\text{(\ref{K7})}}}{=}  \frac{1}{t} \displaystyle{\int^t_0} \omega (s) \ds  \stackrel{{\text{(\ref{K8})}}}{\le} \frac{1}{t} \displaystyle{\int^t_0} (1 + s)^{p - 2} \ds
 = \frac{1}{t} \frac{1}{p - 1}\left\{\left(t + 1\right)^{p - 1} - 1\right\} \, .
\end{eqnarray*}
In case $t \ge 1$ we observe $\{\ldots\} \le (2 t)^{p - 1}$, thus $\frac{1}{t} g' (t) \le c_{53} \  t^{p - 2} \le c_{54} (1 + t)^{p - 2}$, whereas in case $t \le 1$ we use
\[
\frac{1}{t} \left\{\left(t + 1\right)^{p - 1} - 1\right\} = (p - 1) (T + 1)^{p - 2}
\]
for some $T \in (0, 1)$, hence
\[
\frac{1}{t} g' (t) \le c_{55} \le c_{56} (1 + t)^{p - 2}\, ,
\]
and the second inequality in (\ref{K2}) gives the upper bound
\begin{equation}
\label{K9}
D^2 G (\xi) (\eta, \eta) \le c_{57} \left(1 +\left|\xi\right|\right)^{p - 2} |\eta|^2\, .
\end{equation}
With analogous calculations we obtain a lower bound:
\begin{eqnarray*}
\lefteqn{g'' (t) = \omega (t) \ge (1 + t)^{- \mu}, \ t \ge 0\, ,}\\[2ex]
&&\hspace*{-0.75cm}\displaystyle{\frac{1}{t}} g' (t) =  \frac{1}{t} \displaystyle{\int^t_0} \omega (s) \ds  \ge \frac{1}{t} \displaystyle{\int^t_0} (1 + s)^{- \mu} \ds
 = \frac{1}{t(1 - \mu)}\big\{\left(1 + t\right)^{1 - \mu} - 1\big\} \, .
\end{eqnarray*}

\noindent {\bf Case 1:}\ $t \ge t (\mu)$  $(\ge 1$ sufficiently large). Then we have after appropriate choice of $t (\mu)$
\[
\frac{1}{t} g' (t) \ge c_{58} \frac{1}{t} \ge c_{58} (1 + t)^{- \mu }\, .
\]

\noindent {\bf Case 2:}\ $t \le t (\mu)$. Here we observe
\[
\frac{1}{t} \big\{(1 + t)^{1 - \mu} - 1\big\} = (1 - \mu) (1 + \widetilde{T})^{- \mu}
\]
for a suitable $\widetilde{T} \in (0, t (\mu))$, hence
\[
\frac{1}{t} g' (t) \ge (1 + \widetilde{T})^{- \mu} \ge c_{59} \ge c_{59} (1 + t)^{- \mu}\, .
\]
Recalling (\ref{K9}) and (\ref{K2}), the above estimates imply (\ref{G12}) for our density $G$. \qed

\begin{remark}
An equivalent form of (\ref{K8}) is given by
\begin{equation}
\label{K8W}\tag*{$\widetilde{(5.8)}$}
\omega (t) := \Theta (t) (1 + t)^{- \mu}, \ t \ge 0\, , 
\end{equation}
for a continuous function $\Theta : [0, \infty) \to [0, \infty)$ such that $\displaystyle1 \le \Theta (t) \le (1 + t)^{\mu + p - 2}, \ t \ge 0$.

In fact, if $\eta$ is given, let
\[
\Theta (t) := \eta (t) + (1 - \eta (t)) (1 + t)^{p + \mu - 2}\, ,
\]
and if we start from \ref{K8W} we obtain (\ref{K8}) by defining
\[
\eta (t) := \left(\left(1 + t\right)^{p - 2} - \Theta (t) (1 + t)^{- \mu}\right) \left(\left(1 + t\right)^{p - 2} - (1 + t)^{- \mu}\right)^{-1}\, .
\]
\end{remark}
\noindent The density $G$ defined in (\ref{K7}) with $\omega$ as in (\ref {K8}) in general does not satisfy the balancing condition (\ref{K3}): consider $\eta : [0, \infty) \to [0, 1]$ as indicated in the picture below:
\begin{figure}[H]
 \begin{center}
  \includegraphics[scale=1]{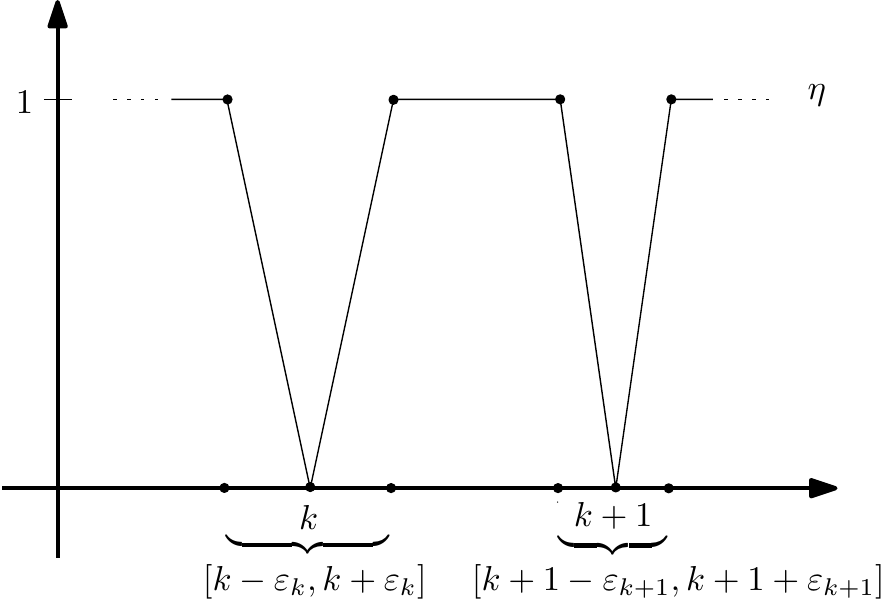}
 \end{center}

\end{figure}

\noindent Here $\varepsilon_k$ denotes a suitable sequence going to zero, we let $\eta \equiv 1$ on each interval $[k + \varepsilon_k, k + 1 - \varepsilon_{k +1}]$ with linear interpolation on $[k - \varepsilon_k, k + \varepsilon_k]$ such that $\eta (k) = 0$ for each $k$. Then it holds for $t \in \mathbb{N}$
\[
t^2 g'' (t) = t^2 (1 + t)^{p - 2}
\]
and at the same time (after appropriate choice of $\varepsilon_k$)
\begin{equation}
\label{K10}
g (s) \le c_{60} \ s
\end{equation}
for $s \ge 0$ sufficiently large, hence (\ref{K4}) is violated. We discuss (\ref{K10}): it holds
\[
g (t) = \int^t_0 \int^s_0 (1 + r)^{- \mu} \dr \ds + \int^t_0 \int^s_0 (1 - \eta (r)) \left\{(1 + r)^{p - 2} - (1 + r)^{- \mu}\right\} \dr \ds\, .
\]
For simplicity let us assume $p \le 2$. Then $\{\ldots\} \le 1$, hence(compare (\ref{G5}) - (\ref{G8}))
\begin{gather*}
g (t) \le c_{61} \left\{t + \int^t_0 \int^s_0 \left(1 - \eta (r)\right) \dr \ds\right\}\le c_{61} \left\{t + \int^t_0 \int^t_0 \left(1 - \eta (r)\right) \dr \ds\right\}\\
\le c _{61} \left\{t + t \int^\infty_0 \left(1 - \eta (r)\right) \dr \right\} \le c_{62} \ t \left\{ 1 + \sum^\infty_{k = 1} \varepsilon_k\right\}\ \, ,
\end{gather*}
and we obtain (\ref{K10}) from the requirement that $\sum^\infty_{k = 1} \varepsilon_k < \infty$. In the case $p > 2$ a slight modification is necessary still leading to (\ref{K10}). In addition, we can choose different functions $\eta$ to get $\lim_{t \to \infty} g(t)/t^q > 0$ for a given number $q \in (1, p)$.
\begin{remark}
As a matter of fact our previous considerations extend to densities
\[
G (\nabla u) = \int^{|\nabla u|}_0 \int^s_0 \left[\eta (r) (\varepsilon_1 + r)^{- \mu} + \left(1 - \eta  (r)\right)(\varepsilon_2 +  r)^{p - 2} \right] \dr \ds
\]
with positive numbers $\varepsilon_1, \varepsilon_2$ and with weight-function $\eta$ as in (\ref{K8}).
\end{remark}
\begin{examp}
We let
\begin{equation}
\label{K11}
G (\xi) := g (|\xi|) := \int^{|\xi|}_0 \int^s_0 (\varepsilon + r)^{\rho (r) - 2} \dr \ds, \ \xi \in \rz\, ,
\end{equation}
with $\varepsilon > 0$ and for a continuous and decreasing function
\begin{equation}
\label{K12}
\rho : [0, \infty) \to [- \mu + 2, p],\quad \rho (0) = p,\quad \lim_{r \to \infty} \rho (r) = 2 - \mu
\end{equation}
with exponents $p, \mu > 1$. Note that (\ref{K11}), (\ref{K12}) can be seen as an approximation of the density $G (\nabla u) = |\nabla u|^{p (|\nabla u|)}$ where $p (|\nabla u|)$ decreases from $p$ to $1$ as $|\nabla u|$ ranges from $0$ to $\infty$, introduced by Blomgren, Chan and Mulet \cite{BCM} for $p = 2$.
\end{examp}
\begin{propos}
The density $G$ from (\ref{K11}) with $\rho$ defined in (\ref{K12}) satisfies the ellipticity condition (\ref{G12}), moreover, the balancing inequality (\ref{K3}) holds.
\end{propos}

\noindent\textit{Proof.} W.l.o.g. we let $\varepsilon = 1$ and observe for any $t \ge 0$ 
\[
 (1 + t)^{- \mu} \le (1 + t)^{\rho (t) - 2} = g'' (t) \le (1 + t)^{p - 2}, 
\]

moreover, it holds
\begin{align*}
&\frac{1}{t}g'(t) = \frac{1}{t} \int^t_0  \big(1 +s\big)^{\rho (s) - 2} \dx\le \frac{1}{t}  \int^t_0  \left(1 + s\right)^{p - 2} \ds \le c_{63} (1 + t)^{p - 2}, \\
&\frac{1}{t}\  g' (t) \ge \frac{1}{t} \int^t_0 (1 + s)^{- \mu} \ds \ge c_{64} (1 + t)^{- \mu},
\end{align*}
we refer to the proof of Proposition 5.4. Thus (\ref{G12}) follows from Proposition 5.1. Next we discuss (\ref{K4}) for $g$ by referring to Proposition 5.3: we have
\begin{eqnarray*}
\lefteqn{g (2t) = \int^{2t}_0  g' (s) \ds = \int^t_0 2 g' (2 s) \ds}\\[2ex]
&&\stackrel{{\text{(\ref{K11})}}}{=} 2 \int^t_0 \int^{2s}_0 (1 + r)^{\rho (r) - 2} \dr \ds = 4 \int^t_0 \int^s_0 (1 + 2r)^{\rho (2r) - 2} \dr \ds\, .
\end{eqnarray*}
In case $p \le 2$ we get by the properties of $\rho$
\begin{align*}
(1 + 2r)^{\rho (2r) - 2} \le (1 + r)^{\rho (2r) - 2},\quad (1 + r)^{\rho (2r) - 2} \le (1 + r)^{\rho (r) - 2}
\end{align*}

(note: $\rho (2r) \le \rho (r)$), hence
\[
g (2 t) \le 4 \int^t_0 \int^s_0 (1 + r)^{\rho (r) - 2} \dr \ds = 4 g (t)\, .
\]
If the value of $p > 1$ is arbitrary, we write
\[
(1 + 2 r) ^{\rho (2 r) - 2} = (1 + r)^{\rho (2 r) - 2} \left(\frac{1 + 2 r}{1 + r}\right)^{\rho (2 r) - 2}
\]
 and use the fact that 
\[
\lim_{r \to \infty} \left(\frac{1 + 2r}{1 + r}\right)^{\rho (2 r) - 2} = 2^{- \mu},
\]
thus $(1 + 2 r) ^{\rho (2 r) - 2} \le c_{65}  (1 + r)^{\rho (2 r) - 2}$ , and by recalling $\rho (2 r) \le \rho (r)$ we obtain as before inequality (\ref{K6}) with a suitable constant. It remains to check (\ref{K5}): we have
\begin{align*}
&g' (t) = {\int^t_0} (1 + s)^{\rho (s) - 2} \ds \ge \int^t_{t/2} (1 + s)^{\rho (s) - 2} \ds\ge {\int^{t}_{t/2}} (1 + s)^{\rho (t) - 2} \ds, 
\end{align*}
since $\rho$ decreases. Writing
\[
\int_{t/2}^t (1 + s)^{\rho (t) -2} \ds = (1 + t)^{\rho (t) - 2} \int_{t/2}^t \left\{\frac{1 + s}{1 + t} \right\}^{\rho (t) - 2} \ds
\]
and observing that
\[
\left\{\frac{1 + s}{1 + t}\right\}^{\rho (t) - 2}  \ge c_{66} > 0 \ \  \mbox{on} \ [t/2, t]\, ,
\]
we see that $g' (t) \ge c_{66} \ \frac{t}{2} (1 + t)^{\rho (t) - 2}$ and (\ref{K5}) is established. \qed


\blfootnote{
\begin{large}\Letter\end{large} Martin Fuchs (fuchs@math.uni-sb.de), Jan Müller (jmueller@math.uni-sb.de)\\
Saarland University (Department of Mathematics), P.O. Box 15 11 50, 66041 Saarbrücken, Germany }
\end{document}